\documentclass[12pt]{article}
\usepackage[english]{babel}
\usepackage{amsmath}
\usepackage[cp1251]{inputenc}
\usepackage{amscd}
\usepackage{amssymb}
\newtheorem{thm}{Theorem}

\newtheorem{cor}{Corollary}
\begin{document}

\begin{center}
%\large
\bf On the precise form of the inverse Markov factor
for convex sets\footnote{{\it MSC} 41A17;
\ {\it key words}: Tur\'{a}n inequality, Markov inequality, convex set}
\end{center}

\begin{center}
{\bf Mikhail A. Komarov}\footnote{Vladimir State University, Gor'kogo street 87,
600000 Vladimir, Russia; \ \ e-mail: kami9@yandex.ru}
\end{center}

{\bf Abstract.}  {\small Let $K\subset \mathbb{C}$ be a
convex compact set, and let $\Pi_n(K)$ be the class of polynomials
of exact degree $n$, all of whose zeros lie in $K$.
The Tur\'an type inverse Markov factor is defined by
$M_n(K)=\inf_{P\in \Pi_n(K)} \left(\|P'\|_{C(K)}/\|P\|_{C(K)}\right)$.
A combination of two well-known results due to Levenberg and Poletsky
(2002) and R\'ev\'esz (2006) provides the lower bound
$M_n(K)\ge c\left(wn/d^2+\sqrt{n}/d\right)$, $c:=0.00015$,
where $d>0$ is the diameter of $K$ and $w\ge 0$ is the minimal width
(the smallest distance between two parallel lines between
which $K$ lies). We prove that this bound is essentially sharp,
namely, $M_n(K)\le 28\left(wn/d^2+\sqrt{n}/d\right)$ for all $n,w,d$.}

\medskip\bigskip

Let $\Pi_n(K)$ be the class of complex polynomials
of exact degree $n$, all of whose zeros lie in a compact
set $K\subset \mathbb{C}$. Let $d=d(K)$ and $w=w(K)$ denote
the diameter of $K$ and the minimal width (i.e., the smallest
distance between two parallel lines between which $K$ lies).
In this note, we investigate the so-called
{\it inverse Markov factor}
\[M_n(K):=\inf_{P\in \Pi_n(K)} \frac{\|P'\|_{K}}{\|P\|_{K}}, \quad n\ge 1
\quad (\|\cdot\|_K:=\|\cdot\|_{C(K)}),\]
for convex compact sets $K$.

Recall that studying the inverse Markov (or Markov-Bernstein)
inequalities was started in 1939 by Tur\'an \cite{Turan}.
In the cases of the disk $K=D:=\{z:|z|\le 1\}$
and the interval $K=I:=[-1,1]$, Tur\'an has established the equality
\[M_n(D)=n/2\]
and the weak equivalence $M_n(I)\asymp \sqrt{n}$,
i.e., the two-sided inequality
\begin{equation}\label{Turan}
    A\sqrt{n}<M_n(I)<B\sqrt{n}
\end{equation}
with some absolute constants $A,B>0$. More precisely,
it was proved that $M_n(I)>\sqrt{n}/6$ ($n\ge 1$)
and $M_n(I)\le \sqrt{n/e}+o(1)$ ($n\to \infty$).

Soon, Er\H{o}d \cite{Erod} has initiated the investigation
of the factor $M_n(K)$ for more general convex sets
(the commented survey of his important paper may be found in \cite[Sec. 1]{Revesz-2013}).
Two remarkable results in this direction were obtained
two decades ago. The first is the result due to
Levenberg and Poletsky, who have proved \cite{Levenberg-Poletsky}
in 2002 the lower bound
\begin{equation}\label{Levenberg-Poletsky}
    M_n(K)\ge \frac{\sqrt{n}}{20d}, \quad d=d(K), \quad n\ge 1,
\end{equation}
where $K$ is {\it an arbitrary convex compact set}. The
order of this bound is sharp and is attained, by (\ref{Turan}),
if $K$ is an interval (a convex compact set of zero width).

The second fundamental result was obtained in 2006 by R\'ev\'esz
\cite{Revesz-2006-JAT}. R\'ev\'esz has proved that if the
width of a convex compact set $K$ is nonzero, $w=w(K)>0$, then
the bound (\ref{Levenberg-Poletsky}) may be significantly
improved as follows:
\begin{equation}\label{Revesz >=}
    M_n(K)\ge 0.0003\frac{wn}{d^2}, \quad n\ge 1.
\end{equation}
He also established the converse inequality, which is valid
for any {\it connected} (not necessarily convex) compact set $K$,
whose width $w$ is nonzero:
\begin{equation}\label{Revesz <=}
    M_n(K)\le 600\frac{wn}{d^2},
    \quad n>2\Big(\frac{d}{16w}\Big)^2\log\frac{d}{16w}
\end{equation}
(in fact, the restriction $n\ge n_0(K):=\max\{1;2(d/16w)^2\log(d/16w)\}$
is assumed).

The estimates (\ref{Revesz >=}), (\ref{Revesz <=}) show the
precise dependence of the quantity $M_n(K)$ on the degree $n$
of polynomials as well as
on the geometric characteristics $w$, $d$ of a convex set $K$
for all {\it sufficiently large} $n$:
\[M_n(K)\asymp wn/d^2.\]
However, in the general case
the formula may not work. Indeed, if we fix $n,d$ and let
$w$ to 0, we obtain $wn/d^2\to 0$, while $K$ 
becomes a segment, $\Delta$ say, of length $d$,
so that $M_n(K)\to M_n(\Delta)>\sqrt{n}/(3d)$ (see (\ref{Turan})).

In this note we fill the gap by proving that R\'ev\'esz's bound
(\ref{Revesz <=}) remains valid without the log-term (Corollary \ref{Cor 1})
and by establishing the precise form of the inverse Markov factor
$M_n(K)$ for arbitrary $n,w,d$ (Theorem \ref{Theorem 1}).

\begin{thm}\label{Theorem 1}
There are absolute constants $c_1,c_2>0$ such that
for any convex compact set $K$ and $n\ge 1$ we have
\begin{equation}\label{main th_0}
  c_1\max\Big\{\frac{wn}{d^2};\frac{\sqrt{n}}{d}\Big\}\le M_n(K)
  \le c_2\max\Big\{\frac{wn}{d^2};\frac{\sqrt{n}}{d}\Big\},
\end{equation}
where $w=w(K)$, $d=d(K)$ and the value $c_2:=28$ is suitable.
The upper bound in $(\ref{main th_0})$ is true for any
connected compact set $K$.
\end{thm}

Obviously, the lower bound in (\ref{main th_0}) immediately
follows from (\ref{Levenberg-Poletsky}), (\ref{Revesz >=}):
\[M_n(K)\ge \max\Big\{0.0003\frac{wn}{d^2};\frac{\sqrt{n}}{20d}\Big\}
\ge c_1 \max\Big\{\frac{wn}{d^2};\frac{\sqrt{n}}{d}\Big\},
\quad \ c_1:=0.0003,\]
and the convexity of $K$ is essential here (see \cite{Levenberg-Poletsky}
for examples of non-convex connected sets whose
inverse Markov factor is $o(\sqrt{n})$ as $n\to \infty$).

Note that the form of both estimates (\ref{Levenberg-Poletsky}),
(\ref{Revesz >=}), as well as of (\ref{main th_0}),
is invariant under affine transformations $t(z)=\alpha z+\beta$
of independent variable, since
$M_n(t(K))=M_n(K)/|\alpha|$, $w(t(K))=|\alpha|w(K)$,
$d(t(K))=|\alpha|d(K)$ (see details in \cite[p. 171]{Revesz-2006-JAT}).

\smallskip

The inequality $(a+b)/2\le \max\{a;b\}\le a+b$ ($a,b\ge 0$)
allows us to write the estimate (\ref{main th_0}) in the following alternative form:
\begin{equation}\label{main th}
  \frac{c_1}{2}\Big(\frac{wn}{d^2}+\frac{\sqrt{n}}{d}\Big)\le M_n(K)
  \le c_2\Big(\frac{wn}{d^2}+\frac{\sqrt{n}}{d}\Big), \quad \ n\ge 1.
\end{equation}
The prototype for the estimate (\ref{main th}) and thus for
Theorem \ref{Theorem 1} is an interesting result
\cite{Erdelyi} obtained by Erd\'elyi around 2005 for the
diamonds $K=S_\varepsilon$ with diagonals
$[-1,1]$ and $[-i\varepsilon,i\varepsilon]$ ($0\le \varepsilon\le 1$)
and for polynomials having certain kinds of symmetry.
In particular, Erd\'elyi has proved the existence of such absolute
constants $c_3,c_4>0$ that
\begin{equation}\label{diamond}
    c_3(n\varepsilon+\sqrt{n})\le
    \inf_{P} \frac{\|P'\|_{S_\varepsilon}}{\|P\|_{S_\varepsilon}}\le
    c_4(n\varepsilon+\sqrt{n}), \quad \ n\ge 1,
\end{equation}
where the infimum is taken over all real-valued polynomials
$P$ from the class $\Pi_n(S_\varepsilon)$. Letting $\varepsilon\to 0$
in (\ref{diamond}), we get the classical Tur\'an's result (\ref{Turan}).
It is not difficult to find that 
$\varepsilon\asymp w(S_\varepsilon)$, namely,
$w(S_\varepsilon)=2\varepsilon/\sqrt{1+\varepsilon^2}$.
For the ellipses with axes $[-1,1]$ and $[-i\varepsilon,i\varepsilon]$,
the result, similar to (\ref{main th}), was presented in \cite{Erod}
without full prove for a lower bound.   

\smallskip

Using Theorem \ref{Theorem 1}, we can significantly improve
the bound (\ref{Revesz <=}). Indeed, if $n>d^2/(kw)^2$ then
$\sqrt{n}/d<kwn/d^2$. Choosing $k=\sqrt{459}$ and using
(\ref{main th_0}) and $28k<600$, we obtain
\[M_n(K)\le 600\frac{wn}{d^2} \quad \text{for all} \ \ n>\frac{d^2}{459w^2},\]
which improves (\ref{Revesz <=}) since $d^2/(459w^2)<n_0(K)$.
In the most interesting case $k=1$ we get

\begin{cor}\label{Cor 1}
For any connected compact set $K$, whose width is nonzero, we have
\[M_n(K)\le 28\frac{wn}{d^2} \quad \text{for} \ \ n>\frac{d^2}{w^2}.\]
\end{cor}

Here the restriction on $n$ is not far from being optimal 
in the sense that if
\[n<d^2/(560w)^2\]
then the bound $M_n(K)\le 28\,wn/d^2$ cannot holds
(for such $n$ we have $28\,wn/d^2<\sqrt{n}/(20d)$,
while $M_n(K)\ge \sqrt{n}/(20d)$ by (\ref{Levenberg-Poletsky})).

More generally, the estimate (\ref{main th_0}) provides a
continuous scale of exponent $1/2\le \gamma\le 1$ in the
formula $M_n(K)\asymp n^\gamma$ depending on the relation
between $n$ and the {\it relative width} $s:=w/d$ of $K$.
For example, $M_n(K)<(28d^{-1})n^{2/3}$ if $s<n^{-1/3}$, and
this bound cannot be improved in the order. It is interesting
that the order of magnitude of the inverse Markov factor for convex domains
may be the same --- $\sqrt{n}$ --- as for a segment.

\begin{cor}\label{Cor 2}
If a natural number $n$ and a convex compact set $K$ satisfy the
condition $w\le d/\sqrt{n}$ then
\[\frac{\sqrt{n}}{20d}\le M_n(K)\le 28\frac{\sqrt{n}}{d}.\]
\end{cor}

{\it Proof of Theorem \ref{Theorem 1}.}
To prove the upper bound in (\ref{main th_0}), we essentially
modify the proof of the estimate (\ref{Revesz <=}) given in \cite{Revesz-2006-JAT}.

Without loss of generality we may assume that a (connected)
compact set $K=K_1$ contains the points $z=+1$, $z=-1$
and its diameter equals $2$:
\[d(K_1)=2, \qquad \pm 1\in K_1.\]
All the points of such a compact set satisfy the simple bounds
\cite[Lemma 1]{Glazyrina-Goryacheva-Revesz}:
\begin{equation}\label{|y|<w}
  |{\rm Im}(z)|\le w(K_1), \quad \ |{\rm Re}(z)|\le 1 \quad \ (z\in K_1).
\end{equation}
For simplicity, we will write $\|\cdot\|$ instead of $\|\cdot\|_{K_1}$.

We need to find polynomials which confirm the correctness of the bound
\[M_n(K_1)\le 28\max\Big\{\frac{wn}{d^2};\frac{\sqrt{n}}{d}\Big\}=7\,\max\big\{wn;2\sqrt{n}\big\},
\quad n\ge 1.\]

If
\begin{equation}\label{Case 1}
    n\le 199 \quad \text{or} \quad w\ge 3/7,
\end{equation}
then the polynomials
\[Q_n(z):=(z-1)^n, \quad Q_n\in \Pi_n(K_1),\]
are suitable. Indeed, let
$\|Q_n'\|=|Q_n'(z_0)|=n|z_0-1|^{n-1}$, $z_0\in K_1$.
We have $|z_0-1|\le d(K_1)=2$. Hence $|Q_n'(z_0)|\le n2^{n-1}$,
while $\|Q_n\|\ge |Q_n(-1)|=2^n$. Therefore,
\[M_n(K_1)\le \|Q_n'\|/\|Q_n\|\le n/2.\]
But if $n\le 199$ or $w\ge 3/7$ then $n=\frac{1}{2}\sqrt{n}\cdot 2\sqrt{n}<8\cdot 2\sqrt{n}$
or $n\le \frac{7}{3}wn$. So,
\[M_n(K_1)<4\,\max\big\{wn;2\sqrt{n}\big\}
\quad \text{in the case (\ref{Case 1})}.\]

Further assume $n>199$ and $w<3/7$.

Let us consider the polynomials
\[p_m(z):=(z^2-1)^m, \quad m\ge 100, \quad p_m\in \Pi_n(K_1),\]
for even $n=2m>199$, and the polynomials
\[P_m(z):=(z-1)(z^2-1)^m, \quad m\ge 100, \quad P_m\in \Pi_n(K_1),\]
for odd $n=2m+1>199$. For all $m$ we have
\begin{equation}\label{||p m||>1}
  \|p_m\|\ge 1 \quad \text{(and} \quad \|P_m\|\ge 1).
\end{equation}
Indeed (see \cite{Revesz-2006-JAT}), since $K_1$ is connected
and both points $\pm 1$ belong to $K_1$, then there is a point
$\tilde{z}\in K_1$ of the form $\tilde{z}=0+i\tilde{y}$. Thus
\[|\tilde{z}-1|\ge 1 \quad \text{and} \quad |\tilde{z}+1|\ge 1,\]
and (\ref{||p m||>1}) follows. We also need another simple estimate
\begin{equation}\label{|p m'|...<2m|z|}
  R(z):=\frac{|p_m'(z)|}{\|p_m\|}\le 2m|z| \quad \text{for every} \ \ z\in K_1.
\end{equation}
In the case $|z^2-1|\ge 1$, this follows immediately:
\[R(z)\le \left|\frac{p_m'(z)}{p_m(z)}\right|=\left|\frac{2mz}{z^2-1}\right|\le 2m|z|,\]
and if $|z^2-1|<1$ then we first use the inequality $\|p_m\|\ge 1$ (see (\ref{||p m||>1})):
\[R(z)\le |p_m'(z)|=2m|z|\,|z^2-1|^{m-1}\le 2m|z|.\]

To obtain Theorem \ref{Theorem 1}, it suffices to prove the following two
assertions concerning the polynomials $p_m$ and the polynomials $P_m$,
respectively.

\bigskip

\noindent{\bf Theorem 1a} {\it \ For any even $n=2m$, $m\ge 100$, we have}
\[\|p_m'\|\le\sqrt{3}\,\max\{wn;2\sqrt{n}\}\,\|p_m\|
\quad \ (w<3/7).\]

\noindent{\bf Theorem 1b} {\it \ For any odd $n=2m+1$, $m\ge 100$, we have}
\[\|P_m'\|\le7\,\max\{wn;2\sqrt{n}\}\,\|P_m\|
\quad \ (w<3/7).\]

{\it Proof of Theorem 1a.} Let us show that
\begin{equation}\label{claim}
  R(z_0)\le \sqrt{3}\,\max\{2m\,w;2\sqrt{2m}\}, \quad z_0=x_0+iy_0,
\end{equation}
for an arbitrary point $z_0\in K_1$ and any $m\ge 100$, where $R(z)$ is defined by (\ref{|p m'|...<2m|z|}).

If $|x_0|\le \sqrt{2}\,w$, then (\ref{claim}) easily
follows by (\ref{|p m'|...<2m|z|}) and $|y_0|\le w$ (see (\ref{|y|<w})):
\[R(z_0)\le 2m|z_0|=2m\sqrt{x_0^2+y_0^2}\le 2m\sqrt{2w^2+w^2}=
\sqrt{3}\cdot 2m\,w,\]
so further we will assume that
\begin{equation}\label{|x_0|>...}
    |x_0|\ge \sqrt{2}\,w.
\end{equation}
Also, instead of (\ref{|p m'|...<2m|z|}), we will apply
another auxiliary estimate (see (\ref{||p m||>1}))
\[R(z_0)=\frac{|p_m'(z_0)|}{\|p_m\|}\le |p_m'(z_0)|=2m\cdot f(x_0,y_0), \qquad z_0=x_0+iy_0,\]
where the function $f$ is defined as
\[f(x,y):=|z(z^2-1)^{m-1}|=\sqrt{x^2+y^2}\,\big(\sqrt{(1+y^2-x^2)^2+4x^2y^2}\big)^{m-1}.\]

Taking (\ref{|x_0|>...}) and (\ref{|y|<w}) into account,
it suffices to estimate the quantity
\[f^*:=\max\{f(x,y): \sqrt{2}\,w\le |x|\le 1, \ |y|\le w\}.\]
The function $f(x,y)$ is even in $x$, is even in $y$, and
increases with $|y|$, therefore
\[f^*=\max\{f(x,w): \sqrt{2}\,w\le x\le 1\}.\]
Put $f(x,w)=\sqrt{g(x^2)}$, where
\[g(t):=(t+w^2)\left(h(t)\right)^{m-1}, \quad
h(t):=(1+w^2-t)^2+4w^2t \quad \ (2w^2\le t\le 1).\]

\medskip

{\it Subcase 1:}
\[1-w\le t\le 1.\]
Here $h(t)\le \max\{h(1-w);h(1)\}$, since the parabola $h(t)$ is convex. But
\[h(1)=4w^2+w^4\le h(1-w)=5w^2-2w^3+w^4\le 5w^2 \quad (w<3/7),\]
so that, for all $1-w\le t\le 1$, we have
\[h(t)\le 5w^2<5\cdot\frac{9}{49}<1, \quad (h(t))^{m-1}\le (h(t))^9
\le 5w^2(45/49)^8\le 2.53\,w^2\]
(by using $m\ge 100>10$) and finally
\[g(t)=(t+w^2)(h(t))^{m-1}\le (1+9/49)\cdot 2.53\,w^2\le 3w^2.\]

\medskip

{\it Subcase 2:}
\[2w^2\le t\le 1-w.\]

First of all, we have the bound
\begin{equation}\label{H(t)<=... (*)}
  h(t)\le 1+2w^2-t, \qquad 2w^2\le t\le 1-w,
\end{equation}
since the convexity of the function $h(t)+t$ gives
\[\max_{t\in [2w^2,1-w]}(h(t)+t)=\max\{h(2w^2)+2w^2; h(1-w)+1-w\},\]
where $h(2w^2)+2w^2=1+9w^4\le 1+2w^2$ and
\[h(1-w)+1-w=(5w^2-2w^3+w^4)+1-w=\]
\[=1+2w^2-w\big((1-w)^3-w^2\big)\le 1+2w^2 \quad \ (0\le w<3/7).\]
Using (\ref{H(t)<=... (*)}), we can estimate $g(t)$ as follows:
\[g(t)\le (t+w^2)(1+2w^2-t)^{m-1}=:g_1(t), \quad \ 2w^2\le t\le 1-w.\]
Direct calculations show that the derivative
\[g_1'(t)=m(1+2w^2-t)^{m-2}(t^*-t), \quad \text{where} \quad t^*:=\frac{1-(m-3)w^2}{m}.\]
We have $t^*\le 1/m\le 1/100$ ($m\ge 100$). Hence
\[t^*<1-w \quad \ (w<3/7).\]
If, moreover, $t^*\le 2w^2$ then $g_1'(t)\le 0$ for all
$t\in [2w^2,1-w]$, and so
\[g(t)\le g_1(2w^2)=3w^2, \qquad 2w^2\le t\le 1-w.\]
If, otherwise, $t^*>2w^2$, then $1+2w^2-t^*<1$ and the maximum
\[\max_{t\in[2w^2,1-w]} g_1(t)=g_1(t^*)=(t^*+w^2)(1+2w^2-t^*)^{m-1}
<t^*+w^2=\frac{1+3w^2}{m}.\]
From this,
\[g(t)<\frac{1}{m}\Big(1+\frac{27}{49}\Big)<\frac{2}{m}, \qquad 2w^2\le t\le 1-w \quad (w<3/7).\]
Note that the inequality $t^*>2w^2$ is equivalent to the condition on $m$:
\[m-1<1/(3w^2).\]

Summarizing the results of Subcases 1 and 2 we get the estimate
\[g(t)\le \max\Big\{3w^2;\frac{2}{m}\Big\}, \qquad 2w^2\le t\le 1,\]
which gives
\[f^*=\max\big\{\sqrt{g(t)}: 2w^2\le t\le 1\big\}\le
\max\Big\{\sqrt{3}\,w;\frac{\sqrt{2}}{\sqrt{m}}\Big\}\]
and
\[R(z_0)\le 2m\cdot f^*\le \sqrt{3}\,\max\Big\{2m\,w;\frac{2\sqrt{2m}}{\sqrt{3}}\Big\}
\le \sqrt{3}\,\max\big\{2m\,w;2\sqrt{2m}\big\}.\]

Thus, we have established (\ref{claim}) over all $z_0\in K_1$.
Applying this to any point at which the norm $\|p_m'\|$ is attained,
we get Theorem 1a.

\medskip

{\it Proof of Theorem 1b.}
Let we take an arbitrary $\alpha>2$ (the exact value of $\alpha$
will be determined later) and an arbitrary point $z\in K_1$.
We have
\[P_m'(z)=(z-1)(z^2-1)^{m-1}(z+1+2mz).\]

If
\[|z-1|\cdot\left|\frac{z+1}{2m}+z\right|\le \alpha|z|\]
then
\[|P_m'(z)|=2m|z-1|\,|z^2-1|^{m-1}\left|\frac{z+1}{2m}+z\right|
\le \alpha\cdot 2m|z|\,|z^2-1|^{m-1}=\alpha|p_m'(z)|.\]
By Theorem 1a,
\[|p_m'(z)|\le\sqrt{3}\,\max\{2mw;2\sqrt{2m}\}\|p_m\|<\sqrt{3}\,\max\{wn;2\sqrt{n}\}\|p_m\|\]
(here $n=2m+1$). But if $\|p_m\|=|p_m(z_0)|$, $z_0\in K_1$,
then the inequalities $|z_0+1|\le d(K_1)=2$ and
$\|p_m\|=\big(|z_0+1||z_0-1|\big)^m\ge 1$ (see (\ref{||p m||>1}))
imply that
\[|z_0-1|\ge 1/2.\]
From this, $\|p_m\|=|z_0^2-1|^m\le 2|z_0-1||z_0^2-1|^m=2|P_m(z_0)|$
and hence
\begin{equation}\label{p m < P m}
    \|p_m\|\le 2\|P_m\|.
\end{equation}
This gives
\[\frac{|P_m'(z)|}{\|P_m\|}\le \frac{2\alpha|p_m'(z)|}{\|p_m\|}
<2\alpha\sqrt{3}\,\max\big\{wn;2\sqrt{n}\big\}.\]

Now let
\[|z-1|\cdot\left|\frac{z+1}{2m}+z\right|>\alpha|z|.\]
Applying the inequality $|z-1|\le d(K_1)=2$, we get
\[\alpha|z|<2\frac{|z+1|}{2m}+2|z|, \quad \text{or} \quad
2m|z|<\frac{2|z+1|}{\alpha-2}.\]
Thus (see (\ref{p m < P m})),
\[|P_m'(z)|\le |z-1|\,|z^2-1|^{m-1}\big(|z+1|+2m|z|\big)
\le |z^2-1|^m\Big(1+\frac{2}{\alpha-2}\Big)\]
\[=|p_m(z)|\frac{\alpha}{\alpha-2}\le \|p_m\|\frac{\alpha}{\alpha-2}
\le 2\|P_m\|\frac{\alpha}{\alpha-2}.\]
But $n\ge 201$ by the assumptions. So, $\sqrt{n}\ge \sqrt{201}>14$ and
\[\frac{|P_m'(z)|}{\|P_m\|}<2\sqrt{n}\cdot\frac{1}{14}\frac{\alpha}{\alpha-2}.\]

To combine the both variants, we choose $\alpha$ by the equality
\[2\alpha\sqrt{3}=\frac{1}{14}\frac{\alpha}{\alpha-2}\]
and put $\alpha=\alpha_0:=2+(28\sqrt{3})^{-1}$. So,
\[\frac{|P_m'(z)|}{\|P_m\|}<2\alpha_0\sqrt{3}\,\max\big\{wn;2\sqrt{n}\big\}
\quad \text{for every} \ \ z\in K_1,\]
where the constant $2\alpha_0\sqrt{3}=4\sqrt{3}+1/14=6.9996\ldots<7$.

Theorem 1b is proved. Hence, Theorem \ref{Theorem 1} is now proved, too.

\end{document}